\documentclass{amsart}

\usepackage{latexsym}

\title[$Q$-Hall-Littlewood functions]{A $q$-analog of Schur's
$Q$-functions }
\author{Geanina Tudose and Michael Zabrocki}
\email{tudose@math.umn.edu, zabrocki@mathstat.yorku.ca }
\address{School of Mathematics, University of Minnesota,
Minneapolis, Minnesota, 55455 and Department of Mathematics and
Statistics, York University, Toronto, Ontario, M3J 1P3}

\font\Sc=cmcsc10

\newtheorem{thm}{Theorem}
\newtheorem{cor}[thm]{Corollary} 

\newtheorem{prop}[thm]{Proposition}
\newtheorem{ex}{Example}
\newtheorem{conj}[thm]{Conjecture}

\begin{document}
\begin{abstract} 
We present a family of analogs of
the Hall-Littlewood symmetric functions in the $Q$-function
algebra. The change of basis coefficients between this family
and Schur's $Q$-functions
 are $q$-analogs of numbers of marked shifted tableaux.
These coefficients exhibit many parallel properties to the Kostka-Foulkes
polynomials.
\end{abstract}
\maketitle 

\def\Z{{\mathbb Z}}
\def\Q{{\mathbb Q}}
\def\C{{\mathbb C}}
\def\N{{\mathbb N}}
\def\D{{\mathcal D}}
\def\horiz{{\mathcal H}}
\def\S{{\bf S}}
\def\H{{\bf H}}
\def\QS{{\bf Q}}
\def\G{{\bf G}}
\def\la{{\lambda}}
\def\smcoeff{{\big|}}
\def\coeff{{\Big|}}
\def\minus{{\epsilon}}
\newdimen\Squaresize \Squaresize=9pt
\newdimen\Thickness \Thickness=0.5pt

\def\Square#1{\hbox{\vrule width \Thickness
     \vbox to \Squaresize{\hrule height \Thickness\vss
        \hbox to \Squaresize{\hss#1\hss}
     \vss\hrule height\Thickness}
\unskip\vrule width \Thickness}
\kern-\Thickness}

\def\Vsquare#1{\vbox{\Square{$#1$}}\kern-\Thickness}
\def\Blk{\omit\hskip\Squaresize}

\def\Young#1{
\vbox{\smallskip\offinterlineskip
\halign{&\Vsquare{##}\cr #1}}}

\newdimen\squaresize \squaresize=5pt
\newdimen\thickness \thickness=0.2pt

\def\square#1{\hbox{\vrule width \thickness
     \vbox to \squaresize{\hrule height \thickness\vss
        \hbox to \squaresize{\hss#1\hss}
     \vss\hrule height\thickness}
\unskip\vrule width \thickness}
\kern-\thickness}

\def\vsquare#1{\vbox{\square{$#1$}}\kern-\thickness}
\def\blk{\omit\hskip\squaresize}
\def\noir{\vrule height\Squaresize width\Squaresize}%

\def\young#1{
\vbox{\smallskip\offinterlineskip
\halign{&\vsquare{##}\cr #1}}}

\def\thisbox#1{\kern-.09ex\fbox{#1}}
\def\downbox#1{\lower1.200em\hbox{#1}}

\section{Introduction}

The space of $Q$-functions, $\Gamma$, is defined to be the
algebra generated by the odd power sum elements $\{p_1, p_3, p_5,
\ldots \}$ as a subalgebra of the space of symmetric functions,
$\Lambda$.
$\Gamma$ is associated to the representation
theory of the spin group and is also related to the projective
represenation theory of the symmetric and alternating groups.
The fundamental basis for this space are Schur's 
$Q$-functions, $Q_\la[X]$,
which are indexed by strict partitions $\la$.
These functions hold the place that the Schur $S$-functions,
$s_\mu[X]$ for
$\mu$ a partition, represent in the algebra of the symmetric
functions.

The space of symmetric functions contains an important basis,
$H_\mu[X;q]$,  the Hall-Littlewood symmetric functions \cite{Hl}.
Through specializations of the parameter $q$ these 
functions interpolate several well studied bases of the symmetric
functions and generalize features of these bases.  They
have elegant properties and may be seen
in many different contexts of combinatorics, algebra,
representation theory, geometry and mathematical physics.

It is natural to ask the question of what the analog of the
Hall-Littlewood symmetric functions in the
$Q$-function algebra should be.
In this paper we introduce a family of functions
$G_\la[X;q] \in \Gamma$ that answers this question since we
observe that this family
shares many of the combinatorial and algebraic properties
of the $H_\mu[X;q]$ functions in the space of symmetric functions.
We expect that these functions will also be interesting from the
perspective of other fields as well. 

From the combinatorial standpoint, we
note that the coefficient of $s_\la[X]$ in the symmetric
function $H_\mu[X;q]$ is the well known Kostka-Foulkes
polynomial.  This family of coefficients are known to be
polynomials in the parameter $q$ with non-negative integer
coefficients and at $q=1$ represent the number of column
strict tableaux of shape $\la$ and content $\mu$.  The
combinatorial tools of jeu de taquin and the plactic monoid
were in part developed to explain the connections between
the Kostka-Foulkes polynomials and the column strict tableaux
\cite{LS}.  

By comparison the coefficient of $Q_\la[X]$ in the function
$G_\mu[X;q]$ is also a polynomial in $q$ and we conjecture
(and prove in certain cases) that it also has coefficients
that are non-negative integers.  At $q=1$ we know that these
coefficients are the number of marked shifted tableaux with
shape $\la$ and content $\mu$.  A version of the RSK-algorithm
was developed by Sagan, Worley and others \cite{H},
\cite{Sa}, \cite{W}, \cite{S1},
 and used to develop the theory of marked shifted
tableaux.  We hope that this theory can be extended to help
answer the question of a combinatorial interpretation for
these coefficients.

Our definition for the functions $G_\la[X;q]$ is motivated by
viewing the symmetric functions $s_\mu[X]$ and $H_\mu[X;q]$
as compositions of operators.  In the case of the Schur
functions, the Bernstein operator $\S_m \in End(\Lambda)$
(\cite{M} p. 96) has the property for $m \geq \mu_1$, 
$$\S_m(
s_\mu[X] ) = s_{(m,\mu_1, \ldots, \mu_{\ell(\mu)})}[X].$$ 
That is, this formula is a recursive defintion for the Schur
functions of degree
$n+m$ as an algebraic relation that raises the degree of a
symmetric function by $m$ acting on a Schur function of degree $n$.
For the Hall-Littlewood symmetric functions, the operator
$\H_m \in End(\Lambda)$ with
$$\H_m( H_\mu[X;q]) = H_{(m,\mu_1, \ldots, \mu_{\ell(\mu)})}[X;q]$$
for $m \geq \mu_1$ is due to Jing \cite{J2}.  In \cite{Z}, it was
noticed that these operators (as well as many others) are related by
a simple algebraic
$q$-twisting, ${\widetilde {\S_m}}^q = \H_m$ (the defintion
of {\,\,}$\widetilde{\,}{\,}^q$ is stated precisely in equation
(\ref{qhat}) below).

The Schur's $Q$-functions, $Q_\la[X]$ may also be seen from this
perspective (\cite{J1}, \cite{LP}, \cite{M} p. 262-3).  That is,
there exists an operator
$\QS_m \in End(\Gamma)$ such that for $m > \la_1$,
$$\QS_m( Q_\la[X] ) = Q_{(m, \la_1, \la_2, \ldots
\la_{\ell(\la)})}[X].$$
Since the defintion of the $q$-twisting {\,\,}$\widetilde{\,}{\,}^q$
extends to $End(\Gamma)$, a natural defintion for an analog
to of the Hall-Littlewood symmetric functions in $\Gamma$ is to
define $\G_m := {\widetilde{\QS_m}}^q$ and then for $m > \la_1$, we
set
$$G_{(m,\la_1, \ldots, \la_{\ell(\la)})}[X;q] :=
\G_m( G_\la[X;q]).$$
This framework provides us only with a possible definition for
the Hall-Littlewood analogs in $\Gamma$.  It remains to show
that these functions share properties similar to those of 
the Hall-Littlewood functions.  In this case we find some
striking similarities that say we have indeed found the
correct analog.

This work is inspired by the results of the
Hall-Littlewood functions and the desire to find analogous structure
in the $Q$-function algebra.  In addition, part of the motivation of
defining these functions and identifying their properties is to find
what features of the Hall-Littlewood symmetric functions are not
unique to the symmetric function algebra and should hold in a more
general setting.  A goal of this reasearch is to possibly identify
what the $q$-twisiting of equation (\ref{qhat}) represents on a
combinatorial, geometric or representation theoretical level
and to show that the $G_\la[X;q]$ are another example of a
structure that seems to exist in a more general context.

The remainder of this paper is divided into three sections and two
appendicies.  The first section is simply a exposition of
definitions and and notation related to the symmetric functions and
$Q$-function algebra.  We develop in some detail the perspective
that bases of the symmetric functions and $Q$-functions can be seen
as compositions of operators that have simple algebraic
definitions.  In the next section we introduce the $G_\la[X;q]$
functions and derive recurrences and some properties that
are analogous to those that exist for the Hall-Littlewood symmetric
functions.  In a final section we discuss a generalization
of the functions $G_\la[X;q]$ that are indexed by a sequence of
strict partitions $(\mu^{(1)}, \mu^{(2)}, \ldots, \mu^{(k)})$ and
the motivation for this generalization.  These functions correspond
to a $q$-analog of the product $Q_{\mu^{(1)}}[X] Q_{\mu^{(2)}}[X]
\cdots Q_{\mu^{(k)}}[X]$ and the coefficients of these functions
correspond to analogs of the generalized or parabolic Kostka
polynomials of \cite{K}, \cite{KS}, \cite{SiWa},  \cite{SW}
and
\cite{SZ}.

Finally, in the first appendix we include tables of transition
coefficients between the $G_\la[X;q]$ basis and the $Q_\la[X]$
basis for degrees $3$ through $9$.  These tables are evidence of a
very strong conjecture that these coefficients are polynomials in
$q$ with non-negative integer coefficients and represent a $q$
analog of the number of marked shifted tableaux.  This suggests
that the marked shifted tableaux should have a poset structure
similar to the charge poset for the column strict tableaux.  In
the second appendix we include a conjectured diagram for a
rank function on the marked shifted tableaux of content $(4,3,2)$. 
This diagram suggests that the structure for the statistic on
marked shifted tableaux is somewhat different than that of the
charge poset even though we conjecture that these statistics
should share many of the same properties.

\section{Notation and Definitions} 

\subsection{Symmetric functions, partitions and columns
strict tableaux} \label{notationsf}

Consider $\Lambda^X$ the ring of series of finite
degree in the variables
$x_1, x_2, x_3, \ldots$ which are invariant under
all permutations of the variables.  This ring is algebraically
generated by the set of elements $\{ p_k[X] = \sum_i x_i^k \}
\subset \Lambda^X$ and hence $\Lambda^X$ is isomorphic to the
ring $\Lambda = \C[p_1, p_2, p_3, \ldots]$ with $deg(p_k) = k$.
We will refer to both $\Lambda$ and $\Lambda^X$ as the
ring of symmetric functions.


$\Lambda$ is a graded ring and basis for the component of
degree $n$ is given by the
monomials $p_\la := p_{\la_1} p_{\la_2} \cdots
p_{\la_{\ell}}$ where $\la$  is a non-increasing
sequence of non-negative integers such that the values
sum to $n$. Such a sequence is called a
partition of $n$ (denoted $\la \vdash n$).
The entries of $\la$ are called the
parts of the partition.
The number of parts that are of size $i$ in $\la$ will
be represented by $m_i(\la)$ and the total number
of non-zero parts
is represented by $\ell(\la) := \sum_i m_i(\la)$ and the
size by $|\la| := \sum_k k m_k(\la) = \sum_i \la_i$. A common
statistic associated to partitions is
$n(\la) := \sum_i (i-1) \la_i$.

The partial order on partitions, $\la\leq\mu$ if and only if
$\sum_{i=1}^k \la_i \leq \sum_{i=1}^k \mu_i$ for all $1 \leq
k \leq \ell(\la)$, is called the dominance order.  We call the operators 
$$
R_{ij} \la = (\la_1, \ldots, \la_i+1, \ldots, \la_j-1, \ldots, \la_{\ell(\la)})
$$
for $1 \leq i \leq j\leq \ell(\la)$ `raising operators' and they have the
property that if $R_{ij} \la$ is a partition, then $R_{ij} \la \geq \la$.

We will consider three additional bases of $\Lambda$ here.
Following the notation of \cite{M}, we define
the homogeneous (complete) symmetric functions as
$h_\la := h_{\la_1}
h_{\la_2} \cdots h_{\la_{\ell(\la)}}$ where
$h_n = \sum_{\la \vdash n} p_\la/z_\la$ and $z_\la =
\prod_{i=1}^{\ell(\la)} i^{m_i(\la)} m_i(\la)!$.
The elementary symmetric functions are $e_\la := e_{\la_1}
e_{\la_2} \cdots e_{\la_{\ell(\la)}}$ where
$e_n = \sum_{\la \vdash n} (-1)^{n-\ell(\la)} p_\la/z_\la$.
By convention we set $p_0 = h_0 = e_0 = 1$ and $p_{-k} =
h_{-k} = e_{-k} = 0$ for $k>0$.
The Schur functions are given by
$s_\la = det\left| h_{\la_i + i - j}\right|_{1 \leq i,j \leq
\ell(\la)}$.  The sets $\{p_\la\}_{\la \vdash n}$, $\{h_\la\}_{\la
\vdash n}$,
$\{e_\la\}_{\la \vdash n}$ and
$\{s_\la\}_{\la\vdash n}$ all form bases for the symmetric functions
of degree $n$.

We will consider the elements of $\Lambda$ as
functors on the space $\Lambda^X$.  If $E$ is an
element of $\Lambda^X$ then we let $p_k[E] = E$ with
$x_i$ replaced with $x_i^k$, so that
$p_k: \Lambda^X \rightarrow \Lambda^X$.  We then
extend this relation algebraically, $p_\la[E]$
will represent the expression $p_{\la_1}[E] p_{\la_2}[E]
\cdots p_{\la_{\ell(\la)}}[E]$.  In particular, if we
take $E = X = x_1 + x_2 + x_3+ \cdots$ then
$p_k[X] = \sum_i x_i^k$ and the map that sends
$p_k \mapsto p_k[X]$ is a ring isomorphism
$\Lambda \rightarrow \Lambda^X$ since
$p_k[X] = x_1^k + x_2^k + x_3^k + \cdots$.

Identities which hold in the ring $\Lambda$ specialize
as well to the ring of symmetric polynomials in a finite
number of variables. We will use the notation $X_n$ to represent
$x_1 + x_2 + \cdots + x_n$, and the map which sends $\Lambda$ to
$\Lambda^{X_n}$ by $p_k \mapsto p_k[X_n]$ corresponds to
the operation of specializing variables in a symmetric series
from an infinite set variables to a symmetric polynomial in a
finite set of variables.  For notational purposes,
capital letters used as variables
the end of the alphabet will be
used to represent a series of variables (e.g. $X = \sum_i x_i$ or
$Y = \sum_i y_i$), while capital letters indexed by a number
will represent a polynomial sum of variables (e.g. $X_n = \sum_{i=1}^n
x_i$ or $Z_n = \sum_{i=1}^n z_i$).

We will need to adjoin to each of the rings $\Lambda, \Lambda^{X}$
and $\Lambda^{X_n}$ a special
element $q$ (or many special elements, if necessary)
which acts much like a variable in this ring, however
$q$ will specialize to values in the field.  $q$
has the special property that $p_k[ q X] = q^k p_k[X]$ and hence
is not an element of our base field since for $c \in \C$, we
have that $p_k[ c X ] = c p_k[X]$.

Notice that by definition
we have in general $p_k[a X + b Y] = a p_k[X] + b p_k[Y]$ for
$a,b \in \C$.  This implies that $f[- X]$ does not represent
the symmetric series $f[X]$ with $x_i$ replaced by $-x_i$ since
$p_k[X] \smcoeff_{x_i \rightarrow -x_i} = (-1)^k p_k[X]$, while
$p_k[-X] =
- p_k[X]$.  To this end we introduce the notation $ f[\epsilon X] =
f[ q X] \smcoeff_{q = -1}$.  In the case of the power sums we have
that $p_k[ \epsilon X ] = (-1)^k p_k[X]$ and hence
$f[ \epsilon X] = f[X] \smcoeff_{x_i \rightarrow -x_i}$.

Consider the series $\Omega = \sum_{n \geq 0} h_n$ which is
not an element of the ring $\Lambda$, but lies in the completion
of this ring.
We will use the morphism $p_k \mapsto p_k[X]$ on this element
as well and manipulations of this notation allow us to derive
the following identities, which we will use repeatedly in our
calculations.
\begin{align}\label{omid1}
\Omega[X] &= \prod_{i} \frac{1}{1-x_i} = \sum_{ n \geq 0 } h_n[X]
\\
\label{omid2}
\Omega[-X]&= \prod_i 1-x_i = 1/\Omega[X]\\
\label{omid3}
\Omega[X+Y] &= \Omega[X]\Omega[Y]\\
\label{omid4}
\Omega[ -\epsilon X] &= \prod_{i} 1+x_i =
\sum_{n\geq0} e_n[X]
\end{align}

Define a generating function of operators $\S(z) = \sum_{m}
\S_m z^m$ where  for an arbitrary symmetric function $P[X]$, $\S(z)
P[X] = P[X - 1/z]
\Omega[zX]$.  In this manner $\S_m$ acts on any symmetric function
raising the degree of the
function by $m$ and has the action
$\S_m P[X] = \S(z) P[X] \coeff_{z^m}$.  A composition of the
operators $\S(z_i)$ produces the expression

\begin{equation} \label{Schurops}
\S(z_1) \S(z_2) \cdots \S(z_k) 1 = \Omega[Z_k X]
\prod_{1 \leq i < j \leq n}
(1-z_j/z_i),
\end{equation}
where $Z_k$ represents the sum $z_1 + z_2+ \cdots +z_k$.  Since
the coefficient of $z^\la$ in $\Omega[Z_k X]$ then it must be
that the coefficient of $z^\la$ in the right hand side of
(\ref{Schurops}) is given by
$ \prod_{1\leq i<j\leq n} (1-R_{ij}) h_{\la}[X],$
where $R_{ij} h_{\la}[X] = h_{R_{ij} \la}[X]$ (considering the
$h$-functions indexed by sequences of numbers).
This is an expression for the Schur function
$s_\la[X]$, hence it follows that $\S_{\la_1} \S_{\la_2} \cdots
\S_{\la_k} 1 = s_\la[X]$.  These operators are due to
Bernstein
\cite{M} p.96.  It follows that
$\S_m(s_\la[X] ) = s_{(m,\la)}[X]$
where $(m,\la)$ denotes
$(m,\la_1, \ldots, \la_{\ell(\la)})$
and
$s_{(a_1, \ldots, a_{\ell})}[X] = det| h_{a_i + i - j}[X] |$.
\vskip7pt
\noindent
{{\bf Remark} 1:}
We follow \cite{M} in the use of raising operators for
our definitions, however we are being imprecise since our
raising operators are not associative or commutative as
defined.  We will consider a symmetric function as a
composition of operators (for example in equation \ref{Schurops})
and the operators $R_{ij}$ serve to raise or lower the
indexing integer of the operator in the $i$ and $j^{th}$ positions
respectively.

By acting on an arbitrary symmetric function using these
operations  and the relations in equations (\ref{omid1})
(\ref{omid2}) and (\ref{omid3}) commutation relations of the
operators follow very nicely.
By expanding the left and right side of the following
expression verifies that
\begin{align*}
\S(z) \S(u) P[X] = - \frac{u}{z}
\S(u) \S(z) P[X].
\end{align*}
By taking the coefficient of $u^m z^n$ in both sides of
the equation, we find that $\S_n \S_m = - \S_{m-1} \S_{n+1}$
which also implies that $\S_m \S_{m+1} = 0$.
Many of the calculations of commutation relations are of
a similar sort of manipulation.

A Young diagram for a partition will be a collection of
cells of the integer grid lying in the first quadrant.
For a partition $\la$, $Y(\la) = \{ (i,j) : 0 \leq j < \ell(\la)
\hbox{ and } 0 \leq i \leq \la_j \}$.  The reason why we
consider empty cells rather than say, points, is because
we wish to consider fillings of these cells.  A tableau is
a map from the set $Y(\la)$ to $\N$, this may be represented
on a Young diagram by writing integers within the cells
of a graphical representation of a Young diagram (see
figure 1).  The shape of the tableau is the partition $\la$.
We say that a tableau $T$ is column strict
if $T(i,j) \leq T(i+1,j)$ and $T(i,j) < T(i,j+1)$ whenever
the points $(i+1,j)$ or $(i,j+1)$ are in $Y(\la)$.  Let
$m_k(T)$ represent the number of points $p$ in $Y(\la)$ such
that $T(p) = k$.  The vector $(m_1(T), m_2(T), \ldots )$
is the content of the tableau $T$.

The Pieri rule describes a combinatorial method for
computing the product of $h_m[X]$ and $s_\mu[X]$
expanded in the Schur basis.
We will use the notation $\la \slash \mu \in \horiz_m$
to represent that $|\la|-|\mu|=m$ and for
$1 \leq i \leq \ell(\la)$, $\mu_i \leq \la_i$ and
$\mu_i \geq \la_{i+1}$.
It may be easily shown that
\begin{equation}
h_m[X] s_\mu[X] = \sum_{\la\slash\mu \in
\horiz_m} s_\la[X].
\end{equation}

This gives a method for computing the expansion of the
$h_\mu[X]$ basis in terms of the Schur functions.
Consider the coefficients $K_{\la\mu}$ defined by
the expression
\begin{equation}\label{Kostka}
h_\mu[X] = \sum_{\la \vdash |\mu|} K_{\la\mu} s_\la[X].
\end{equation}
$K_{\la\mu}$ are called the Kostka numbers
and are equal to the number of column strict
tableaux of shape $\la$
and content $\mu$.

\subsection{Kostka polynomials and Hall-Littlewood
symmetric functions} Define the following symmetric functions
\begin{align}
H_\la[X;q] &= \prod_{i<j} \frac{1-R_{ij}}{1-q R_{ij}} h_\la[X]
\nonumber\\ &=\label{eq:defHla}
\prod_{1\leq i<j\leq n} (1 +(q-1)R_{ij}+ (q^2-q) R_{ij}^2 + \cdots)
h_{\la}[X].
\end{align}
They will be referred to as Hall-Littlewood symmetric functions
as they are transformations of the symmetric polynomials
defined by Hall \cite{Hl}  (see  \cite{M} for a modern account
where $Q_\mu(x;q)$ in their notation is $H_\mu[X(1-q);q]$ in
ours).  The coefficient of
$s_\la[X]$ in
$H_\mu[X;q]$ is known as the Koskta Foulkes polynomial
$K_{\la\mu}(q)$.  That is, we have the expansion
\begin{equation}
H_\mu[X;q] = \sum_\la K_{\la\mu}(q) s_\la[X].
\end{equation}
We will present some of the properties of the Kostka-Foulkes
polynomials and the Hall-Littlewood symmetric functions below.
First, it will be important to establish some identities for
manipulating these functions.

Let $\H(z) = \sum_m \H_m z^m$ be defined as the operation
$\H(z) P[X] = P[ X - (1-q)/z] \Omega[z X]$.  Taking the
coefficient of $\H_m = \H(z) \coeff_{z^m}$ has the effect of
raising the degree of the symmetric function it is acting on
by $m$.  A composition of these operators has the expression
\begin{equation} \label{eq:compHm}
\H(z_1) \H(z_2) \cdots \H(z_k) 1 =
\Omega[Z_k X] \prod_{1\leq i < j\leq k} \frac{1-z_j/z_i}{1-q
z_j/z_i}.
\end{equation}
Now since $h_{R_{ij}\la} = \frac{z_j}{z_i} \Omega[Z_k X]
\coeff_{z^\la}$, it is clear that the coefficient of
$z^\la$ in the
right hand side of (\ref{eq:compHm}) is exactly the right
hand side of (\ref{eq:defHla}) and hence  $\H_{\la_1}
\H_{\la_2} \cdots \H_{\la_{\ell(\la)}} 1 = H_\la[X;q]$. This
operator also satisfies the relations
$\H_{m-1} \H_m = q \H_m \H_{m-1}$ and $\H_{m-1} \H_n-
q \H_m \H_{n-1} = q \H_n \H_{m-1} -  \H_{n-1} \H_{m}$.  This
relation
can be derived as we did for the Schur function operators
by demonstrating $(z - q u) \H(z) \H(u) = (q z - u) \H(u)
\H(z) $ on an arbitrary symmetric function $P[X]$.

This family of operators $\H_m$ is due to Jing \cite{J2}
and they are sometimes referred to as `vertex operators'
for the Hall-Littlewood symmetric functions.

For an element $V \in Hom(\Lambda, \Lambda)$, define
\begin{equation} \label{qhat}
{\widetilde V}^q P[X] = V^Y P[q X + (1 -q) Y] \coeff_{Y=X},
\end{equation}
where $V^Y$ denotes that as an operation on symmetric
functions in the $Y$ variables only and $Y=X$ represents
setting the $Y$ variables equal to the $X$ variables after
the operation is completed.
This is a $q$-analog of the operator $V$ and we remark
that ${\widetilde {\S(z)}}^q = \H(z)$.  This follows by calculating
\begin{align*}
  {\widetilde {\S(z)}}^q P[X] &= \S^Y(z) P[ q X + (1-q) Y]
\coeff_{Y=X}\\ &=  P[q X + (1-q) (Y - 1/z) ] \Omega[ zY]
\coeff_{Y=X}\\ &= P[ X - (1-q)/z] \Omega[ z X] = \H(z) P[X].
\end{align*}
This relationship between $\S(z)$ and $\H(z)$
is the motivation for our definition of the $q$-analog of
Schur's $Q$-functions.

The functions $H_\la[X;q]$ interpolate between the functions
$s_\la[X] = H_\la[X;0]$ and $h_\la[X] = H_\la[X;1]$.
The Kostka-Foulkes polynomials are defined as the
$q$-polynomial coefficient of
$s_\la[X]$ in $H_\mu[X;q]$ and hence we have the expansion
analogous to (\ref{Kostka}).
\begin{equation}
H_\mu[X;q] = \sum_{\la \vdash |\mu|} K_{\la\mu}(q) s_\la[X].
\end{equation}
The coefficients $K_{\la\mu}(q)$ are clearly polynomials
in $q$, but it is surprising to find that the coefficients of
the polynomials are non-negative integers.

A defining
recurrence can be derived $K_{\la\mu}(q)$ in terms of
the Kostka-Foulkes polynomials indexed by partitions of size
$|\mu| - \mu_1$ using the formula for $\H_m$.
This recurrence is often referred to as the `Morris
recurrence' for the Kostka-Foulkes polynomials \cite{Mo}.  
The action of $\H_m$ on the Schur functions is given by
\begin{equation} \label{Hmaction}
\H_m(s_\mu[X]) = \sum_{i\geq0} \sum_{\mu\slash\la \in
\horiz_i} q^i s_{(m+i,\la)}[X].
\end{equation}
It is not immediately obvious that at $q=1$, the previous
equation reduces to the Pieri rule and at $q=0$ the formula
is simply $\S_m(s_\mu[X]) = s_{(m,\mu)}[X]$.  Using
(\ref{Hmaction}) and equating coefficients of $s_\la$ on both
sides of the equation $\H_m(H_\mu[X;q]) = \sum_\la
K_{\la\mu}(q) \H_m(s_\la[X])$, we arrive at the Morris
recurrence

\begin{equation}\label{morrisKlamu}
K_{\alpha, (m,\mu)}(q)= \sum_{s=1}^{t: \alpha_t \geq m} (-1)^{s-1}
q^{\alpha_s-m} \sum_{ \la: \la/\alpha^{(s)} \in
     \mathcal{H}_{(\alpha_s-m)} }
       K_{\la \mu}(q),
\end{equation}
where $m > \mu_1$ and $\alpha^{(s)}$ is $\alpha$ with part
$\alpha_s$ removed.

The Kostka-Foulkes
polynomials and the
generating functions $H_\mu[X;q]$ have the following
important properties which we simply list here so that
we may draw a connection to analogous formulae.  For a
more detailed reference of these sorts of properties
we refer the interested reader to the excellent survey
article \cite{DLT}.

\begin{enumerate}

\item\label{prop3} the degree in $q$ of $K_{\la\mu}(q)$ is
$n(\mu)-n(\la)$.

\item\label{prop4} $K_{\la\mu}(0) = \delta_{\la\mu}$ which implies
$H_\mu[X;0] = s_\mu[X]$, $K_{\la\mu}(1) = K_{\la\mu},$ so
that
$H_\mu[X;1] = h_\mu[X]$, $K_{\la\la}(q) = 1$ and
   $K_{(|\mu|)\mu}(q) = q^{n(\mu)}$.  We also have that $K_{\la\mu}(q)=
0$ if $\la<\mu$.

\item\label{prop2} $K_{\la\mu}(q) = \sum_T q^{c(T)}$, where
the sum is  over
all column strict tableaux of shape $\la$ and content $\mu$
and $c(T)$ denotes the charge of a tableau $T$
(see \cite{LS}) and hence is a polynomial with non-negative
integer coefficients.

\item A combinatorial interpretation for these coefficients
exists in terms of objects known as rigged configurations
\cite{KR}.

\item\label{prop8} $H_{(1^n)}[X;q] = e_n\left[ \frac{X}{1-q} \right]
(q;q)_n$ where $(q;q)_n = \prod_{i=1}^n (1-q^i)$.

\item\label{prop9} If $\zeta$ is $k^{th}$ root of unity,
$H_{\mu}[X;\zeta]$ factors into a product of symmetric functions.

\item\label{prop10} Set $K_{\mu \lambda}'(q):=  q^{n(\lambda)-n(\mu)}
K_{\mu \lambda}(1/q)$, then
$
K_{\mu \lambda}'(q) \geq K_{\mu \nu}'(q)
   \mbox{  for } \lambda \leq \nu .
$

\item\label{prop11} $K_{\la+(a),\mu+(a)}(q) \geq
K_{\la,\mu}(q)$, where $\la+(a)$ represents the partition
$\la$ with a part of size $a$ inserted into it.

\item\label{prop12} $K_{\la\mu}(q)
= \sum_{w \in S_n} sign(w) {\mathcal P}_q( w(\la +
\rho) - (\mu+\rho))$ where ${\mathcal P}_q( \alpha)$ is the
coefficient of $x^\alpha$ in $\prod_{1 \leq i < j \leq n}
(1 - q x_i/x_j)^{-1}$,
a $q$ analog of the Kostant partition function and
$\rho = (\ell(\mu) -1, \ell(\mu)-2, \ldots, 1, 0)$.

\item\label{prop14} $H_\mu[X;q] H_\la[X;q] = \sum_{\gamma}
d_{\la\mu}^\nu(q) H_\nu[X;q],$ for some coefficients
$d_{\la\mu}^{\nu}(q)$ with the property that
if the Littlewood-Richardson coefficient $c_{\la\mu}^{\nu}
= 0$ then $d_{\la\mu}^\nu(q) = 0$.  These coefficients are
a transformation of the Hall algebra structure coefficients.

\item\label{prop15} For the scalar product
$\left< s_\la[X], s_\mu[X] \right>
= \delta_{\la\mu}$, the basis $H_\mu[X(1-q);q]$ is orthogonal with
respect to $H_\la[X;q]$, that is $\left<
H_\la[X;q], H_\mu[X(1-q);q] \right> = 0$ if $\la \neq \mu$.

\end{enumerate}

\begin{figure}
$$\Young{4&6&7\cr3&5&5\cr2&2&3&4&6\cr1&1&1&1&2&3\cr}
\hskip .5in
\Young{\Blk&\Blk&\Blk&5\cr\Blk&\Blk&3'&3&4&5\cr
\Blk&2'&2&3'&4'&4\cr1'&1&2'&2&2&3&3\cr}
$$
\begin{caption}
{The diagram on the left represents a
column strict tableau of shape $(6,5,3,3)$
and content $(4,3,3,3,2,2,2,1)$.  The diagram on the right
represents a shifted marked tableau of shape $(7,5,4,1)$ and
content $(4,4,4,3,2)$.  This tableau has labels which are
marked on the diagonal. }
\end{caption}
\end{figure}

\subsection{Schur's $Q$-functions, strict partitions,
and marked shifted tableaux}

The $Q$-function algebra is a sub-algebra
of the symmetric functions $\Gamma =
\C[p_1, p_3, p_5, \ldots]$. A typical monomial in
this algebra will be $p_\la$, where $\la$ is a partition
and $\la_i$ is odd.
A partition $\la$
is strict if $\la_i > \la_{i+1}$ for all $1
\leq i \leq \ell(\la)-1$ and a partition $\la$ is odd if
$\la_i$ is odd for $1 \leq i \leq \ell(\la)$.
We will use the
notation $\la \vdash_s n$ (respectively $\la \vdash_o n$)
to denote that $\la$ is a partition of size $n$ that is
strict (respectively odd).
Note that the number of
strict partitions of size $n$ and the number of odd partitions of
size $n$ is the same (proof: write out a generating function for each
sequence).

The analog of the homogeneous and elementary symmetric
functions in $\Gamma$ are the functions
$q_\la := q_{\la_1} q_{\la_2} \cdots
q_{\la_{\ell(\la)}},$
where $q_n =
\sum_{\la \vdash_o n} 2^{\ell(\la)} p_\la/z_\la$.
Define an algebra morphism $\theta : \Lambda \rightarrow
\Gamma$ by the action on the $p_n$ generators as
$\theta(p_n) = (1-(-1)^n) p_n$.  That is $\theta(p_n) = 2 p_n$ if
$n$ is odd and $\theta(p_n) = 0$ for $n$ even.  $\theta$
has the property that $\theta( h_n) = \theta(e_n) = q_n$
and may be represented in our notation as
$\theta(p_n[X]) = p_n[(1-\minus)X]$.
Under this morphism, our Cauchy element may also be
considered a generating function for the $q_n$ elements
since
\begin{equation}
\Omega[(1-\minus)X] = \sum_{n\geq0} q_n[X] =
\prod_i \frac{1+x_i}{1-x_i}.
\end{equation}

It follows that $\{ p_\la \}_{\la \vdash_o n}$,
$\{ q_\la \}_{\la \vdash_o n}$,
$\{ q_\la \}_{\la\vdash_s n}$ are all bases for the
subspace of $Q$-functions of degree $n$.  Another
fundamental basis for this space are the Schur's
$Q$-functions $Q_\la[X] = \theta(H_\la[X;-1])$.  These
functions hold a similar place in the $Q$-function algebra
that the Schur functions hold in $\Lambda$. In particular,
$\{ Q_\la[X] \}_{\la\vdash_s n}$ is a basis for the
$Q$-functions of degree $n$.

In analogy with the Schur functions,
$Q_\la[X]$ may also be defined with a raising operator
formula by setting $q=-1$ and applying the $\theta$
homomorphism to equation (\ref{eq:defHla}). We
arrive at the formula:

\begin{equation}
Q_\la[X] = \prod_{i<j} \frac{1-R_{ij}}{1+ R_{ij}} q_\la[X]
= \prod_{i<j} (1 - 2 R_{ij} + 2 R_{ij}^2 - \cdots)
q_\la[X], \label{ropGq}
\end{equation}
where the operators now act as $R_{ij} q_\la[X]
= q_{R_{ij} \la}[X]$.  Furthermore, they have a formula
as the coefficient in a generating function:

\begin{equation}
Q_\la[X] = \Omega[(1-\minus)Z_n X]
\prod_{1 \leq i < j \leq n} \frac{1 - z_j/z_i}{1 + z_j/z_i}
\coeff_{z^\la}.
\end{equation}

As with Schur functions and
the Hall-Littlewood functions, the raising operator
formula leads us to an operator definition.
By setting $\QS(z) P[X]
= P\left[X - \frac{1}{z} \right]\Omega[(1-\minus)z X]$,
it is easily shown that
$$\QS(z_1) \QS(z_2) \cdots \QS(z_n) 1 = \Omega[(1-\minus)Z_n X]
\prod_{1 \leq i < j \leq n} \frac{1 - z_j/z_i}{1 + z_j/z_i},
$$ and hence if we set $\QS_m P[X] = \QS(z) P[X] \coeff_{z^m}$
then $\QS_m( Q_\la[X] ) =
Q_{(m,\la)}[X]$ as long as $m > \la_1$.  The commutation
relations for the $\QS_m$ are
\begin{equation}\label{m-n}
\QS_m \QS_n = - \QS_n \QS_m \hbox{ for }m\neq -n,
\end{equation}
\begin{equation}\label{m-m}
\QS_m \QS_{-m} = 2(-1)^m - \QS_{-m} \QS_m \hbox{ if }m \neq 0,
\end{equation}
\begin{equation}\label{0}
\QS_m^2 = 0\hbox{ if }m \neq 0\hbox{ and }
\QS_0^2 = 1.
\end{equation}
These formulas allow us to straighten the $Q_\mu[X]$
functions when they are not indexed by a strict
partition.

The $Q$-function algebra is endowed with a natural
scalar product.   If we set $\left< p_\la, p_\mu \right>_\Gamma =
2^{\ell(\la)} \delta_{\la\mu} z_\la$ for $\la, \mu \vdash_o n$,
then it may be shown that we also have
\begin{equation}
\left< Q_\la[X], Q_\mu[X] \right>_\Gamma  = 2^{\ell(\la)}
\delta_{\la\mu}.
\end{equation}

A shifted Young diagram for a partition will again be
a collection
of cells lying in the first quadrant.  For a strict
partition $\la$, let $YS(\la) = \{(i,j):0 \leq j
\leq \ell(\la)\hbox{ and }j-1 \leq i \leq \la_j + j -1\}$.
A marked shifted tableau $T$ of shape $\la$ is a map from
$YS(\la)$ to the set of marked integers $\{ 1' < 1 < 2' <
2 < \ldots \}$ that satisfy the following conditions
\begin{itemize}
\item $T(i,j) \leq T(i+1,j)$ and $T(i,j) \leq T(i,j+1)$

\item If $T(i,j) = k$
for some integer $k$ (i.e. has an unmarked label)
then $T(i,j+1) \neq k$

\item If $T(i,j) = k'$ for some marked label $k'$ then
$T(i+1,j) \neq k'$.
\end{itemize}

    We may represent these objects graphically
with a Young diagram representing $\la$ and the cells filled
with the marked integer alphabet.  If $T$ is a marked
shifted tableau, then we will set $m_i(T)$ as the number
of occurrences of $i$ and $i'$ in $T$.  The sequence
$(m_1(T), m_2(T), m_3(T), \ldots )$ is the content of
$T$.

The combinatorial definition of the marked shifted
tableaux is defined so that it reflects the change of
basis coefficients between the $q_\mu$ and $Q_\la$
basis.  The rule for computing the product of
$q_m[X]$ and $Q_\mu[X]$ when expanded in the Schur
$Q$-functions is the analog of the Pieri rule for
the $\Gamma$ space.
If $\la \slash \mu \in \horiz_m$ then $a(\la\slash\mu)$
will represent $1 +$ the number of $1 < j \leq \ell(\la)$ such that
$\la_j>\mu_j$ and $\mu_{j-1}>\la_j$.
We may show that
\begin{equation}\label{Qpieri}
q_m[X] Q_\mu[X] = \sum_{\la\slash\mu \in \horiz_m}
2^{a(\la\slash\mu) - \ell(\la) + \ell(\mu)}
Q_\la[X].
\end{equation}
Denote by $L_{\la\mu}$ the number of marked shifted tableaux $T$
of shape $\la$ and content $\mu$ (where $\la$ is a
strict partition) such that $T(i,i)$ is not a marked integer.
We may expand
the function $q_\mu[X]$ in terms of the $Q$-functions
using (\ref{Qpieri}) to show
\begin{equation}\label{L}
q_\mu[X] = \sum_{\la\vdash|\mu|} L_{\la\mu} Q_\la[X].
\end{equation}


\section{The $Q$-Hall-Littlewood basis $G_\lambda(x;q)$ for the algebra
$\Gamma$}
\noindent
In this section we  define a new family of functions $G_\la[X,q]$
which seems to play the same role as the Hall-Littlewood functions
$H_\la[X,q]$  in the $Q$-functions algebra. These functions are
introduced via a raising operator formula similar
to~(\ref{eq:defHla}). This definition permits an equivalent
interpretation via a corresponding vertex operator $\G_m$ whose
properties are analogues to both the Hall-Littlewood vertex operator
$\H_m$ and $\QS_m$.

Note: From here, unless otherwise stated, all partitions are considered
strict.

\subsection{Raising operator formula}

We define the following analog of the Hall-Littlewood functions in
the subalgebra $\Gamma$
\begin{equation}\label{def}
G_\la [X;q] :=\prod_{1\leq i<j\leq n}
\left(\frac{1+qR_{ij}}{1-
   qR_{ij}}\right) \left(\frac{1-R_{ij}}{1+ R_{ij}}\right)
   q_\la[X]=\prod_{1\leq i<j\leq n}
\left(\frac{1+qR_{ij}}{1-qR_{ij}}\right) Q_\la[X].
\end{equation}
We call the functions
  $G_\la \in \Gamma \otimes_{\mathbb{C}}
\mathbb{C}(q) $  the {\it $Q$-Hall-Littlewood  functions}.

In $\Gamma \otimes \mathbb{C}(q)$ this family can be
expressed in the basis of $Q$-functions as
\begin{equation}\label{eqfund}
G_\mu [X;q] =\sum_{\la} L_{\la \mu}(q) Q_\la[X],
\end{equation}
\noindent
which can be viewed as a $q$-analog of~(\ref{L}).
We call the coefficients $L_{\la \mu}(q)$ the {\it $Q$-Kostka 
polynomials}. We shall see
  that this family of polynomials shares many of the same properties with
  the classical Kostka-Foulkes polynomials. Tables of these
  coefficients are given in an Appendix. It follows from~(\ref{def}) that
  $L_{\la \mu}(q)$ have integer
coefficients and
$L_{\la\mu}(q)=0$ if $\la<\mu$.
This shows
\begin{prop}
The $ G_\lambda$, $\lambda$ strict, form a $\mathbb{Z}$-basis for
$\Gamma \otimes_{\mathbb{Z} } \mathbb{Z}(q)$.
\end{prop}
  \noindent
The basis $G_\la$ interpolates between the Schur's
$Q$-functions and the functions $q_\mu$
because $ G_\la [X;0]= Q_\la[X]$
and $ G_\la[X;1]= q_\la[X]$ as is clear from~(\ref{def}).

Since the coefficient of $z^\la$ in $\Omega[(1-\minus)Z_n
X]$ is $q_\la[X]$ equation (\ref{def}) implies
\begin{equation}\label{GFG}
G_\la [X;q] = \prod_{1\leq i<j\leq n} \left(
\frac{1-z_j/z_i}{1+z_j/z_i}\right)\left(
\frac{1+qz_j/z_i}{1-qz_j/z_i}\right)
\Omega[(1-\epsilon)Z_n X]
\coeff_{z^\la}.
\end{equation}
Define  the operator $\G(z)$  acting on an arbitrary symmetric
function $P[X]$ as
\begin{equation}\label{defop}
\G(z)P[X]=
P\left[X-\frac{1-q}{z}\right]\Omega[(1-\epsilon)zX].
\end{equation}
The operator $\G(z)$ defines a family of operators as $\G(z) = \sum_{m \in \Z}
\G_m z^m$ and hence $\G_m P[X] = \G(z) P[X] \coeff_{z^m}$.

If we consider a composition of these operators acting on
the symmetric function $1$ we obtain
the following
\begin{equation}
\G(z_1) \G(z_2) \cdots \G(z_n) 1 =
\prod _{1\leq i<j\leq n} \left(
\frac{1-z_j/z_i}{1+z_j/z_i}\right)\left(
\frac{1+qz_j/z_i}{1-qz_j/z_i}\right)
\Omega[(1-\epsilon)Z_nX],
\end{equation}
which together with relation~(\ref{GFG}) gives
$$
G_\lambda[X;q]= \G_{\la_1}\ldots \G_{\ell(\la)}(1).
$$
\vskip11pt
Next we investigate some properties of this operator. First, $\G_m$
satisfies the following commutation relation.

\begin{prop}\label{com}
For all $r,s \in \mathbb{Z}$ we have
$$
(1-q^2)(\G_r \G_s + \G_s \G_r)+ q(\G_{r-1}\G_{s+1} - \G_{s+1}\G_{r-1}
      + \G_{s-1}\G_{r+1} - \G_{r+1}\G_{s-1})
$$
$$
  = 2(-1)^r(1-q)^2
      \delta_{r,-s}.
$$
\end{prop}

\noindent
{\bf Proof}
We will prove this relation in a few steps.
Consider $\G(u)$ and $\G(z)$ the operator $\G$ defined above on the 
variable $u$ and  $z$ respectively. We are looking at the composition of 
these two operators.

{\bf Step 1.} We may write
$$
\G(u)\G(z)= G(z,u)F(z/u)F(-qz/u),
$$
where $G(z,u)$ is an operator  symmetric in $z$ and $u$ defined by
$$
G(z,u)P[X]:=P\left[ X-\frac{1-q}{u}-\frac{1-q}{z}\right]
\Omega[(1-\epsilon)(z+u)X]$$
  and $\displaystyle{F(t):= \frac{1-t}{1+t}}$. This is easily seen
from
\begin{equation}\label{later}
\G(u)\G(z)P[X]= \G(u)P\left[ X-\frac{1-q}{z}\right]
\Omega[(1-\epsilon)z X]=
\end{equation}
$$
= P\left[X-(1-q)\left(\frac{1}{z}+\frac{1}{u}\right)\right]
\Omega[(1-\epsilon)(z+u)X]
\Omega\left[(1-\epsilon)(q-1)\frac{z}{u}\right].
$$
Note that first two factors are exactly $G(u,z)$  and $
\displaystyle{
\Omega\left[(1-\epsilon)(q-1)\frac{z}{u}\right]  }$ is equal to
$F(z/u)F(-qz/u)$.  Denote by $\alpha(t):= F(t)+F(t^{-1})
= \sum_{r \in \Z} 2 (-1)^r t^r$.

\vskip 6pt
{\bf Step 2.}
Let $q_i^\perp$ be the adjoint of multiplication by $q_i[X]$ in
the algebra $\Gamma$ under the scalar product $<\,,\,>_\Gamma$.
Consider $Q^\perp(t)$ its generating function, i.e., $Q^\perp(t) =
\sum_iq_i^\perp t^i$.

It is not difficult to see that $P[X+t]= Q^\perp(t)P[X]$. In fact it
suffices to show this  for a suitable basis element in $\Gamma$,
namely the power sum $p_\la[X]$, where $\la$ is an odd partition. We
show that $q_i^\perp p_\la[X]= p_\la[X+t]\coeff_{t^i}$.  We have
\begin{equation}\label{plethqip}
p_\mu[X+t]\coeff_{t^i}= (p_{\mu_1}[X]+t^{\mu_1})\ldots
(p_{\mu_\ell}[X]+t^{\mu_\ell})\coeff_{t^i}.
\end{equation}
At the same time we also know that $p_\la^\perp p_\mu =
2^{-\ell(\la)}\frac{ z_\mu}{
z_\nu} p_\nu$ if $\nu \uplus \la = \mu$
since
\begin{equation}
p_\la^\perp p_\mu\coeff_{p_\nu} = \left< p_\la^\perp p_\mu,
\frac{2^{\ell(\nu)} p_\nu}{z_\nu }
\right>_\Gamma = \frac{2^{\ell(\nu)}}{z_\nu }
\left<  p_\mu,p_{\la\uplus
\nu} \right>_\Gamma = \frac{2^{\ell(\nu)} z_\mu}{ 2^{\ell(\mu)}
z_\nu} \delta_{\mu,\la\uplus\nu}.
\end{equation}
Now because $q_i = \sum_{\la \vdash_o i} 2^{\ell(\la)} p_\la/ z_\la$,
then $q_i^\perp p_\mu[X] = \sum_{\la \uplus \nu = \mu}
\frac{z_\mu}{z_\la z_\nu} p_\nu$ where $\la \vdash_o i$ and
this is exactly the right hand side of (\ref{plethqip}).

Also recall that
$\Omega[(1-\epsilon)tX]= Q(t)$, where $Q$ is the generating functions 
for $q_i$. If we replace these in the expression of $G(z,u)$ we obtain
$$
G(z,u)=
Q(z)Q(u)Q^{\perp}(-u^{-1})Q^{\perp}(-z^{-1})
Q^{\perp}(qu^{-1})Q^{\perp}(qz^{-1}).
$$

\noindent
We know that $Q(u)Q(z)\alpha(u/z)= \alpha(u/z)$ (\cite{M} Chap
III. 8, p. 263), and hence we have as well
$Q^\perp(u)Q^\perp(z)\alpha(u/z)=
\alpha(u/z)$. Therefore
  $$G(z,u)\alpha(u/z)=\alpha(u/z).$$

\vskip6pt
{\bf Step 3.}
Consider the following expression.
\begin{align*}
&[\G(u)\G(z)(1-qz/u)(1+qu/z)+ \G(z)\G(u)(1-qu/z)(1+qz/u)]=\\
=& (1+qz/u)(1+qu/z)G(u,z)\alpha(u/z)\\
=& (1+qz/u+qu/z+q^2)\alpha(u/z),
\end{align*}
and so
\begin{align*}
\G(u)\G(z)&(1-qz/u+qu/z-q^2)+ \G(u)\G(z)(1-qu/z+qz/u-q^2) \\
&= (1+qz/u+qu/z+q^2)\alpha(u/z) = (1-q)^2\alpha(u/z).
\end{align*}
since we have as a series $ t \alpha(t) = - \alpha(t)$.
By taking the coefficient of $u^rz^s$ in the left hand
side when expanded as series and equating with the corresponding
coefficient in the right hand side we obtain the desired relation.

\hfill
$\Box$

For $q=0$ in the relation above we recover the
commutation relations of the operator $\QS$ given in equations
(\ref{m-n}), (\ref{m-m}) and (\ref{0}).   At $q=1$, $\G_m$ becomes
multiplication by $q_m$ and hence is commutative.

\vskip 11pt
Formula (\ref{GFG}) may be used to derive the action of the
operator $\G_m$ on the basis of Schur's
$Q$-functions.

\begin{prop}\label{pactq} For $m>0$,
\begin{equation}\label{actq}
\G_m (Q_\lambda[X])= \sum_{i \geq 0} q^i \sum_{\mu: \lambda/\mu \in
     \mathcal{H}_i} 2^{a(\lambda/\mu)} (-1)^{\epsilon(m+i,\mu)}
     Q_{\mu+(m+i)}[X],
\end{equation}
where $\mu +(k)$ denotes the partition formed by adding a
part of size $k$ to
    the partition $\mu$, and $\epsilon(k, \mu)+1$ represents which part
    $k$ becomes in $\mu+(k)$ ($Q_{\mu+(k)}[X]=0$ if $\mu$ contains
a part of size $k$).  For
$m\leq0$ a similar statement can be made using the commutation
relations (\ref{m-n}), (\ref{m-m}) and (\ref{0}).
\end{prop}

\noindent
{\bf Proof}
From~(\ref{GFG}) the action of $\G_m$ on a function $P[X] \in \Gamma$
can be  written as
\begin{align*}
\G_m P[X]=& P[X-(1-q)/z]\Omega[(1-\epsilon)zX]\coeff_{z^m}\\
=& \sum _{i \geq 0} q^iz^{-i}(q_i^\perp P)[X-1/z]
\Omega[(1-\epsilon)zX]\coeff_{z^m}\\
=& \sum _{i \geq 0} q^i(q_i^\perp P)[X-1/z]
\Omega[(1-\epsilon)zX]\coeff_{z^{m+i}},
\end{align*}
since $P[X+t]=\sum_{i \geq 0}q_i^\perp P[X]t^i$. Thus
$$
\G_m Q_\la[X]= \sum_{i \geq 0} q^i \QS_{m+i} (q_i^\perp Q_\la[X]),
$$
\noindent
where $q_i^\perp$ applied to $Q_\la$ is
$$
q_i^\perp Q_\la[X] = \sum_{\mu: \lambda/\mu \in \mathcal{H}_i}
2^{a(\lambda/\mu)} Q_\mu[X],
$$
If $m > 0$, equation~(\ref{actq}) follows from (\ref{m-n}) and
(\ref{0}). If $m \leq 0$ we may need to use the commutation
relation (\ref{m-m}) for $\QS_n$ to straighten the index to a strict
partition. In general we have $\QS_{m+i}(Q_\mu[X])=
Q_{(m+i,\mu)}[X]$, where we prepend the part $(m+i)$ (possibly
negative) to the partition
$\mu$.
\hfill
$\Box$
\vskip11pt

\begin{ex}
We compute $G_{(3,2,1)}[X;q]$ using the proposition above. We have
\begin{align*}
G_{(3,2,1)}[X;q]=& \G_3 (\G_2(Q_{(1)}[X])) \\
=& \G_3\left( \sum_{i \geq 0} \sum_{(1)/\mu \in \mathcal{H}_i}
2^{a((1)/\mu)} (-1)^{\epsilon(2+i,\mu)} Q_{\mu +(2+i)} [X] \right )
\\
=& \G_3( Q_{(2,1)}[X]) + 2q\G_3( Q_{(3)}[X] )\end{align*}
\begin{align*}
  =& \sum_{i \geq 0} 
\sum_{(2,1)/\mu
   \in \mathcal{H}_i}
2^{a((2,1)/\mu)} (-1)^{\epsilon(3+i,\mu)} Q_{\mu +(3+i)} [X] +\\
  &\qquad + 2q \left ( \sum_{i \geq 0} \sum_{(3)/\nu
   \in \mathcal{H}_i}
2^{a((3)/\nu)} (-1)^{\epsilon(3+i,\nu)} 
Q_{\nu +(3+i)} [X] \right )\\
=& (q^0 2^0 Q_{(3,2,1)} + q^1 2^1 Q_{(4,2)}+ q^2 2^1
Q_{(5,1)})\\ &\quad +2q (q^1  2^1
Q_{(4,2)} + q^2 2^1 Q_{(5,1)} + q^3 2^1 Q_{(6)})\\
=& Q_{(3,2,1)} + (2q +4q^2) Q_{(4,2)}+ (2q^2+ 4q^3)Q_{(5,1)}
+ 4q^4 Q_{(6)}.
\end{align*}
\end{ex}

\subsection {Properties of the polynomials $L_{\la \mu}(q)$}

The $Q$-Kostka polynomials introduced here have a number of
remarkable properties that are very similar to those of
Kostka-Foulkes polynomials listed in the previous section.
We have already seen the analog of
Property~\ref{prop4} holds for $Q$-Kostka
polynomials. In what follows we will consider some of the
other remaining properties.

An important consequence of equation~(\ref{actq}) is a
Morris-like recurrence
which expresses the Q-Kostka polynomials $L_{\la \mu}(q)$ in terms of
smaller ones.

\begin{prop}
We have the following recurrence
\begin{equation}\label{morris}
L_{\alpha, (n,\mu)}(q)= \sum_{s=1}^{t: \alpha_t \geq n} (-1)^{s-1}
q^{\alpha_s-n} \sum_{ \la: \la/\alpha^{(s)} \in
     \mathcal{H}_{(\alpha_s-n)} }2^{a(\la/\alpha^{(s)}) }
       L_{\la \mu}(q),
\end{equation}
where $n > \mu_1$ and $\alpha^{(s)}$ is $\alpha$ with part
$\alpha_s$ removed.
\end{prop}

\noindent
{\bf Proof} If $n > \mu_1$ we have that
   \begin{equation}\label{first}
\G_n G_\mu[X;q] = G_{(n,\mu)}[X;q] = \sum_\alpha L_{\alpha, 
(n,\mu)}(q)Q_\alpha[X].
\end{equation}
On the other hand $G_\mu[X;q] = \sum_\la L_{\la \mu}(q)Q_\la[X]$ and so
$$
\G_n \left ( \sum_\la  L_{\la \mu }(q)Q_\la[X] \right )=
\sum_\mu L_{\la \mu}(q) \G_n (Q_\la[X]).
$$
Using the action in~(\ref{actq}) we have
\begin{equation}\label{latter}
\G_n G_\mu[X;q]= \sum_\la L_{\la \mu}(q)
\sum_{i \geq 0} q^i \sum_{\nu: \la/\nu \in
     \mathcal{H}_i} 2^{a(\la/\nu)} (-1)^{\epsilon(n+i,\nu)}
     Q_{\nu+(n+i)}[X].
\end{equation}
   For $\alpha= \nu +(n+i)$, equating the coefficients of $Q_\alpha$
   in~(\ref{first}) and (\ref{latter})  we
   get
$$
L_{\alpha, (n,\mu)}(q)=\sum_\la \sum_{i \geq 0} q^i 
2^{a(\la/\alpha-(n+i))}
(-1)^{\epsilon(n+i,\alpha-(n+i))} L_{\la \mu}(q).
$$
By reindexing $i:=\alpha_s-n$ for $\alpha_s-n \geq 0$ we
obtain the desired recurrence~(\ref{morris}).
\hfill
$\Box$
\vskip5pt
\begin{ex}
Let $n=5$ and $L_{(6,2), (5,2,1)}(q) = 2q +4q^2$. Using the recurrence
we have one $s$ such that $\alpha_s \geq 5$, i.e. $\alpha_1=6$. So
$$
L_{(6,2), (5,2,1)}(q)= q^{6-5} \sum_{ \la/(2) \in \mathcal{H}_1 } 2^
{ a(\la/(2))} L_{\la
   (2,1)}(q)$$
$$
= q \big{(}  2L_{(21), (21)}(q) + 2 L_{(3),(21)}(q) \big{)}= q( 2 + 
2\cdot
2q) = 2 q + 4 q^2.$$
\end{ex}

\vskip10pt
As a consequence of the Morris-like recurrence we have the following

\begin{cor}\label{cor2}
   Let $\mu \leq \la$ in dominance order.

1. If $n >\la_1$
then $L_{(n, \la), (n,\mu)}(q)=L_{\la\mu}(q)$.

2. $L_{\lambda \lambda}(q)=1$ and  $L_{(|\lambda|)
\lambda}(q)= 2^{\ell(\la)-1}q^{n(\lambda)}$.

3. $2^{\ell(\mu)-\ell(\la)}$ divides $L_{\la \mu}(q)$.
\end{cor}

\noindent
{\bf Proof} 1. There is only one term in the recurrence~(\ref{morris})
in this
     case which is exactly $L_{\la \mu}(q)$.

2. The first is a consequence of (1). For the second, we have that the 
only term on the right hand side is
     $q^{|\lambda|-\lambda_1} 2 L_{(|\lambda|-\lambda_1)
     (\lambda_2,\ldots)}(q)$ which by induction is
     $q^{|\lambda|-\lambda_1 + n((\lambda_2, \ldots))}2\cdot 
2^{\ell(\la)-2} =
    2^{\ell(\la)-1} q^{n(\lambda)}$.
This is the analog of Property~\ref{prop4} for the Kostka-Foulkes 
polynomials.

3. We use induction on $\ell(\mu)$ and  the Morris recurrence to derive 
this property.

\noindent
If $\ell(\mu)=1$ we know that $G_{(m)}[X;q]=Q_{(m)}[X]$ and
thus $\la$ can  be only $(m)$
  and the assertion holds.

For the induction step we use equation~(\ref{morris}). We need to show 
that $2^{\ell(\mu)-\ell(\alpha)+1}$ divides $L_{\alpha, (n,\mu)}(q)$. 
 From the induction hypothesis we know
that  $2^{\ell(\mu)-\ell(\la)}$ divides $L_{\la \mu}(q)$ for every such 
$\la$
  on the left hand side of equation~(\ref{morris}).

  The partitions $\la $  have $\ell(\la) \in\{ \ell(\alpha), 
\,\ell(\alpha)-1\}$. If $\ell(\la) = \ell(\alpha)-1$ we are done.
If $\ell(\la) = \ell(\alpha) $  then $\la \neq \alpha^{(s)}$, since 
$\alpha^{(s)}$ has length one less than $\alpha$,  and thus 
$a(\la/\alpha^{(s)}) \geq 1$. This implies that  
$2^{\ell(\mu)-\ell(\la)+1}$ divides each $2^{a(\la/\alpha^{(s)})}L_{\la, 
\mu}(q)$ and thus $L_{\alpha,(n,\mu)}(q)$.
\hfill
$\Box$
\vskip8pt
Using the Morris recurrence we are also able to obtain a formula for
the degree  of $L_{\la\mu}(q)$ similar to
Property~\ref{prop3} for the Kostka-Foulkes polynomials.
\begin{prop}
If $\mu \leq \la$ in dominance order, we have
$$
\hbox{deg}_q L_{\la \mu}(q) = n(\mu) - n(\la).
$$
\end{prop}
\noindent
{\bf Proof}
   We prove this assertion by induction on $\ell(\mu)$.
For $\ell(\mu)=1$, the equality is obvious.

For the induction step we use the recurrence~(\ref{morris}).
Fix now an index $s$ on the right hand side of the 
equation~(\ref{morris}).
  Denote by $\mu_{(0)}^{(s)}= (\alpha_1+\alpha_s-n, 
\alpha_2,\ldots,\alpha_{s-1},\alpha_{s+1}, \ldots)$.

We claim that $n(\mu_{(0)}^{(s)}) <n(\la)$ for any other $\la$ such that 
$\la/\alpha^{(s)} \in \mathcal{H}_{(\alpha_s-n)} $. We have that
$$
n(\mu_{(0)}^{(s)})=\sum_{i=2}^{s-1}(i-1) \alpha_i
+\sum_{i \geq 0}(s+i-1)\alpha_{s+i+1}
$$
while
$$
n(\la)= \sum_{i=2}^{s-1}  (i-1)(\alpha_i+\epsilon_i)
+\sum_{i \geq 0}(s+i-1)(\alpha_{s+i+1}+\epsilon_{s+i}),
$$
where $\la= (\alpha_1+\epsilon_1,\ldots,\alpha_{s-1}+\epsilon_{s-1},
\alpha_{s+1}+\epsilon_s, \ldots)$ and $\sum_i
\epsilon_i=\alpha_s-n$. Moreover if $\la \neq\mu_{(0)}^{(s)} $
there exists at least one $\epsilon_i$ with $i\geq 2$ such
that $\epsilon_i\neq 0$.

Therefore  $n(\mu_{(0)}^{(s)}) <n(\la)$. We thus have proved that among 
the polynomials $L_{\la \mu}(q)$ in the second sum of~(\ref{morris}), 
the polynomial
$L_{ \mu_{(0)}^{(s)}, \mu}(q)$ has the highest  degree, namely 
$n(\mu)-n(\mu_{(0)}^{(s)})$.

Next we show that in the first sum, the highest degree is obtained
for  $s=1$. That is to say $deg_q
\left(q^{\alpha_1-n}L_{\mu_{(0)}^{(1)}, \mu} (q) \right) > deg_q
\left(q^{\alpha_s-n}L_{\mu_{(0)}^{(s)}, \mu}(q)\right)$, hence
$$
\alpha_1-n+n(\mu)-n(\mu_{(0)}^{(1)}) > \alpha_s- n +n(\mu)
-n(\mu_{(0)}^{(s)}),$$
$$
\alpha_1+\sum_{i=2}^{s-1}(i-1)\alpha_i +\sum_{j \geq s+1}(j-2)\alpha_j >
\alpha_s + \sum_{i=3}^{s-1}(i-2)\alpha_i +\sum_{j \geq s}(j-2)\alpha_j
$$
which is
$$
\alpha_1 + \alpha_2 + \sum_{i=3}^{s-1} (i-1)\alpha_i >
\sum_{i=3}^{s-1} (i-2)\alpha_i+(s-1)\alpha_s=
$$
$$
=(\alpha_3+\alpha_s) +(2\alpha_4+\alpha_s)+\cdots
+[(s-3)\alpha_{s-1}+\alpha_s] +2\alpha_s.
$$
The last inequality  is true as $(i-1)\alpha_i > (i-2)\alpha_i+s$
for $i=3,\ldots, s-1$ and  $\alpha_1+\alpha_2 >2\alpha_s$.

Thus we have $deg_q L_{\alpha, (n,\mu)}(q) =
(\alpha_1-n)+n(\mu)-n(\mu_{(0)}^{(1)})$. Finally we need to show
this is in fact $n((n,\mu))-n(\alpha)$. That is
$$
  \alpha_1-n +\sum _{i \geq 2}(i-1)\mu_i - \sum_{i \geq 3}(i-2)\alpha_i =
\sum_{i\geq 0} \mu_i -\sum_{i \geq 2} (i-1) \alpha_i,
$$
and by simplifying we obtain $\sum \alpha_i-n =\sum \mu_i$, which
is obviously true.

Hence $deg_q \left( L_{\alpha, (n,\mu))}(q) \right)=
n((n,\mu))-n(\alpha) $ and  the proof is complete.
\hfill
$\Box$

\vskip10pt
The property that is most suggestive that these polynomials are
analogs of the Kostka-Foulkes polynomials is

\begin{conj}\label{conj1}
The Q-Kostka polynomials $L_{\la \mu}(q)$ have non-negative
coefficients.
\end{conj}

We will prove this conjecture for some particular cases.
In general we believe that there should exist a similar combinatorial
interpretation as for the Kostka-Foulkes polynomials. More precisely
there should exist a statistic function $d$ on the set of marked shifted
tableaux, similar to the charge function on column strict
tableaux, such that
$$
L_{\la \mu}(q)= \sum_T q^{d(T)}
$$
summed over marked shifted tableaux of shifted shape $\la$ and content
$\mu$ with diagonal entries unmarked.

In addition, we conjecture that this function must have the property 
that if
$T$ and $S$ are two marked shifted tableaux such that by erasing the 
marks the
two resulting tableaux coincide, then $d(T)= d(S)$.

For some
of the polynomials $L_{\la\mu}(q)$, this observation
determines completely the statistic on the tableaux.  For
instance there are two marked shifted tableaux classes of
shape $(5,3,1)$ and content $(4,3,2)$ and
$L_{(5,3,1),(4,3,2)}(q) = 2q + 4q^2$.  Clearly the tableau
with a $3$ in the first row must have statistic $1$ and
with $3$ in the second row has statistic $2$.  On the other
hand, $L_{(8,1),(4,3,2)}(q) = 4 q^5 + 4 q^6$.  This
polynomial does not uniquely determine which of the two
tableaux have statistic $5$ and $6$.  We have used the
function $G_{(4,3,2)}[X;q]$ to draw a conjectured tableau
poset (similar to the case of column strict tableau) for
the marked shifted tableaux with unmarked diagonals of
content $(4,3,2)$ in an appendix.

Another intriguing property of this statistic function $d$ is that
the values it takes are not too different than the charge function.
It seems that in general we have that  for given $\la$ and $\mu$ the
set of $\{ d(T), T\hbox{ in the summation of $L_{\la \mu}(q)$} \}$
is a subset of $\{ c(T), \hbox {$T$ column strict tableaux of shape
$\la$ and content $\mu$} \}$, where $c$ is the usual charge.  This
suggests that there should be relationship between these two
statistics; however, we have so far not been able to establish what
that link might be.

\begin{prop}
(1) For $\mu$ a two-row partition and  $\la > \mu$ we have
$L_{\la\mu}(q) = 2 q^{n(\mu) - n(\la)}$.

(2) If $\mu$ has the property $\displaystyle{\mu_i \geq \sum_{j
\geq i+1} \mu_j}$, Conjecture~\ref{conj1} is true.
\end{prop}

\noindent
{\bf Proof}
(1). Let us consider
$G_\mu[X;q]$ and let
$\mu=(n,m)$. We have that $G_{(n,m)}[X;q]=\G_n(Q_{(m)}[X]) $
which by Proposition~\ref{pactq} is
$$
\G_n (Q_{(m)}[X])= \sum_{i \geq 0} q^i \sum_{\mu: (m)/\mu \in
     \mathcal{H}_i} 2^{a((m)/\mu)} (-1)^{\epsilon(n+i,\mu)}
     Q_{\mu+(n+i)}[X] .
$$
 From this we deduce that $i=0,1,\ldots, m$ and $\mu=(m-i)$. Thus
$$
\G_n (Q_{(m)}[X])= Q_{(n,m)}[X] +\sum_{i=1}^m 2 q^i Q_{(n+i, m-i)}[X]
$$
and the proof is complete.

(2). In this case we prove it by induction on $\ell(\mu)$ and using
the Morris recurrence~(\ref{morris}).

The case $\ell(\mu)=1$ is clear as $L_{\la\mu}(q)= \delta_{\la
\mu}$. For the induction step  consider $L_{\alpha (n,\mu)}(q)$ as
in the right hand side of~(\ref{morris}).  Under our assumption
there is just one index in the first sum, i.e., only $\alpha_1$ can
be greater than $n$. This is true since $|\alpha|=n+|\mu|,\  \alpha
\geq (n,\mu)$ in dominance order and $n \geq
  |\mu|$. Thus the right hand side does not contain negative signs
and by induction it is non-negative . Hence $L_{\alpha (n,\mu)} (q)$
has non-negative coefficients.
\hfill
$\Box$

\vskip10pt
We also  note that monotonicity properties, similar to
Property~\ref{prop10} and~\ref{prop11}, hold for the $Q$-Kostka
polynomials.
\begin{conj}\label{mon}
Let $L_{\la \mu}'(q):=  q^{n(\mu)-n(\la)} L_{\la
     \mu}(q^{-1})$. We have
$$
L_{\la \mu}'(q) \geq  2^{\ell(\nu)-\ell(\mu)} L_{\la \nu}'(q),
\qquad \mbox{for } \mu
\leq \nu \hbox{ in dominance order}.
$$
\end{conj}
We can prove this fact by using induction and the
recurrence~(\ref{morris}) for the case $\mu_1= \nu_1$.
\begin{ex}
Let $\la = (6,2), \ \mu = (4,3,1), \ \nu = (5, 2,1 )$. We have $n(\la)=
2, \ n(\mu)= 5$, and $n(\nu)= 4$. The  $L'$ polynomials are
$$
  L'_{\la \mu}= q^{5-2}( 4/q^2 + 4/q^3) = 4+4q, \qquad L'_{\la \nu}=
  q^{4-2} (2/q+4/q^2)= 4+2q,
$$
and thus $L'_{\la \mu}(q) \geq 2^{3-3} L'_{\la \nu}(q)$.
\end{ex}

\vskip 10pt
Another property of the Kostka-Foulkes polynomials case that seems
to hold in our case refers to the growth of the polynomials $L$. For
the Kostka-Foulkes polynomials the conjecture belongs to Gupta
(see \cite{DLT} and references therein).

\begin{conj}
If $r$ is an integer that is not a part in either partitions
$\lambda$ or $\mu$, then
$$
L_{\la + (r), \mu+(r)}(q)\geq L_{\la \mu }(q).
$$
\end{conj}
The case where $r>\la_1$ (which also ensures that $r > \mu_1$) is
obviously true since $L_{(r, \la), (r,\mu)}(q)=L_{\la
    \mu}(q)$ (see Corollary~\ref{cor2}).

\begin{ex}
Let $\la = (5,3), \mu = (4,3,1) $ and $a= 2$. We have
$$
  L_{(5,3,2), (4,3,2,1)}(q)-L_{(5,3), (4,3,1)}(q)= 2q +4q^2+8q^3 -
(2q+4q^2) = 8q^3.
$$
\end{ex}

\vskip7pt
\subsection{Another expression for $L_{\la\mu}(q)$.}The polynomials
$L_{\la\mu}(q)$ have a similar interpretation to
property~\ref{prop12} using an analog of the
$q$-Kostant partition function. We follow the construction in~\cite{DLT}.
In order to write equation~(\ref{ropGq}) as
\begin{equation}\label{q}
q_\la [X]= \prod_{i<j} \Big{(}\frac{1+R_{ij}}{1- R_{ij}}\Big{)}^{-1}
Q_\la[X].
\end{equation}
we will use linear maps from the group algebra $\mathbb{Z}[\mathbb{Z}^n]$
to the algebra $\Gamma$. A basis of $\Z[\Z^n]$ will consist of formal
exponentials
$\{e^\alpha\}_{\alpha \in }$ which satisfy relations $e^\alpha
e^\beta=e^{\alpha+\beta}$. In fact we identify the ring with the
ring of Laurent polynomials in $x_1,\ldots, x_n$ and set $e^\alpha=
x^\alpha$. With this in mind we are viewing all our polynomials in
$\Gamma_n$ (or $\Lambda_n$) as linear homomorphisms from
$\Z[\Z^n]$ to $\Gamma$ i.e.
$$
Q:e^\lambda \rightarrow Q(e^\lambda)= Q_\la \quad q:e^\lambda
\rightarrow q(e^\lambda)= q_\la.
$$
If we now set $\displaystyle{ \zeta_n := \prod_{1\leq i<j\leq n}
   \left ( \frac{1+x_i/x_j}{1-x_i/x_j}\right ) } $,
we have that $\displaystyle{\zeta_n = \sum_{\alpha \in \mathbb{Z}^n}
\mathcal{R}(\alpha)e^\alpha }$
where $\mathcal{R}(\alpha) = \sum_t a_t 2^t$ and $a_t$
counts the number of ways the vector
$\alpha$ can be written as a sum of positive roots of type
$A_{n-1}$, $t$ of which are
distinct. The positive roots in the root lattice of $A_{n-1}$ are
$\{e_i-e_j\}_{1\leq i < j \leq n}$, where $e_i= (0,\ldots,1,\ldots
0)$ is the canonical basis of $\Z^n$.

Since $q(e^\la)=Q(\zeta_n e^\la)$ we have that
$$
q_\la[X]=q(e^\la)= Q(\sum_{\alpha \in
\mathbb{Z}^n}\mathcal{R}(\alpha)e^\alpha e^\la)
$$
$$
= \sum_{\alpha \in \mathbb{Z}^n}\mathcal{R}(\alpha)Q_{\la+\alpha}[X].
$$
\noindent
  If we consider the same argument for  $\displaystyle
{ G_\la[X,q]=\prod_{1\leq i<j\leq n}
\left(\frac{1+qR_{ij}}{1-qR_{ij}}\right) Q_\la[X]} $ we need to define 
the  $q$-analog of $\zeta_n$ as
$$
\zeta_n(q):= \prod_{i<j} \left ( \frac{1+qx_i/x_j}{1-qx_i/x_j}\right ),
$$
and thus $
\displaystyle{ \zeta_n (q)= \sum_{\alpha \in \mathbb{Z}^n}
\mathcal{R}_q(\alpha)e^\alpha }$, where $\mathcal{R}_q(\alpha)=
\sum_{t,k} a_{t,k} 2^t q^k$ and  $a_{t,k}$ counts the number of ways
the vector
$\alpha$ can be written as a sum of $k$ positive roots, $t$ of
which are distinct.   Hence
$$
G_\la[X,q]= \sum_{\alpha \in \mathbb{Z}^n}
\mathcal{R}_q(\alpha)Q_{\la+\alpha}[X].
$$
This yields another expression for the $Q$-Kostka polynomials in
terms of
$\mathcal{R}_q(\alpha)$ as
$$
L_{\lambda\mu} (q)=\sum_{\alpha :\,Q_{\alpha+\mu}=\pm 2^r Q\lambda}
\pm 2^r\mathcal{R}_q(\alpha).
$$
The index of the sum reflects the straightening of a $Q$-function
indexed by an integer sequence and it is a consequence of the
commutation relations~(\ref{m-m}). It is possible to express the
equation above using the action of the symmetric group on Schur's
$Q$-functions, yielding an alternating sum  similar to
Property~\ref{prop12}. Unfortunately the action of the symmetric
group on Schur's
$Q$-functions indexed by a general integer vector
is not as elegant
as for Schur functions.

\vskip7pt
\noindent
{\bf Remark:}~~
Most of the properties of the $Q$-Kostka polynomials
$L_{\la\mu}(q)$ are analogous to the Kostka-Foulkes polynomials,
but a few properties do not seem to generalize.
\begin{enumerate}


\item The analog of Property~\ref{prop9} does not seem to
hold since computations of $G_\la[X;q]$ where $q$ is set to
a root of unity do not factor.

\item There does not seem to exist an elegant
relationship between $G_\la[X;q]$ and its dual basis
(Property~\ref{prop15}).

\item A property similar to that of Property~\ref{prop14} does
not seem to hold.  We do not know if there is a relationship between
$G_\la[X;q]$ and a Hall-like algebra.

\item
We do not know if an analog of the Macdonald
symmetric functions should exist.   A
family of functions which mimic the formulas for the
Macdonald symmetric functions in \cite{Z} may easily be defined,
but the specializations of the variables indicate that the
same sort of
positivity and symmetry properties of the coefficients cannot
hold through this definition.

\end{enumerate}

\section{Generalized (parabolic) $Q$-Kostka polynomials}

There exists in the literature a few generalizations of the
Kostka-Foulkes polynomials that correspond to $q$-analogs of
multiplicities of irreducibles in tensor products of
irreducible representations
(Littlewood-Richardson coefficients).  In
\cite{SZ} formulas were introduced for realizing
`generalized' or `parabolic' Kostka coefficients \cite{SW}
as coefficients
appearing in families of symmetric functions defined
as compositions of operators.  This
construction may also be extended to the $Q$-function
algebra providing a means of defining a generalization
of the $Q$-Kostka polynomials that corresponds to a
$q$-analog of coefficients in products of $Q_\mu[X]$.

Let $\mu^* = (\mu^{(1)}, \mu^{(2)}, \ldots, \mu^{(k)})$
be a sequence of partitions and ${\bar {\mu^*}}$ the
concatenation of all those partitions. Take $\eta =
(\ell(\mu^{(1)}), \ell(\mu^{(2)}), \ldots,
\ell(\mu^{(k)}))$ and $n = \sum_i \eta_i$, then set
$Roots_\eta =
\{(i,j) : 1 \leq i \leq \eta_1 + \cdots + \eta_r <
j \leq n \hbox{ for some } r\}$.

Now set
\begin{equation}
H_{\mu^*}[X;q] = \prod_{(i,j) \in Roots_\eta}
\frac{1}{1-q R_{ij}} s_{\bar {\mu^*}}[X].
\end{equation}

The parabolic Kostka polynomials are defined as
the coefficients of the Schur basis
in these symmetric functions. That is, the polynomials
$K_{\la;\mu^*}(q)$ are defined by the coefficients in
the expansion
\begin{equation}
H_{\mu^*}[X;q] = \sum_{\la \vdash |\mu^*|}
K_{\la;\mu^*}(q) s_\la[X].
\end{equation}

The functions $H_{\mu^*}[X;q]$ and the parabolic
Kostka
coefficients have the following properties.

\begin{itemize}
\item If ${\bar {\mu^*}}$ is a partition then it
is conjectured that $K_{\la;\mu^*}(q)$ has non-negative
integer coefficients (in certain cases this is known).

\item $H_{\mu^*}[X;0] = s_{\bar {\mu^*}}[X]$.
${\bar {\mu^*}}$ need not be a partition, but this
is consistent with the definition of $s_\la$ in section
\ref{notationsf}.

\item $H_{\mu^*}[X;1] = s_{\mu^{(1)}}[X] s_{ \mu^{(2)}}[X]
   \cdots s_{\mu^{(k)}}[X]$ and in this sense the coefficients
$K_{\la;\mu^*}(q)$ are $q$-analogs of the Littlewood-Richardson
coefficients

\item If $\mu^* = ((\gamma_1),
(\gamma_2),
\ldots, (\gamma_{\ell(\gamma)}))$ where $\gamma$ is a
partition, then $H_{\mu^*}[X;q] = H_\gamma[X;q]$.

\item If $\bar {\mu^*}$ is a partition then
$H_{\mu^*}[X;q] = s_{\bar {\mu^*}}[X] + \sum_{\la > {\bar
{\mu^*}}} K_{\la;\mu^*}(q) s_\la[X]$.

\item There exists an operator $\H_\gamma$ such that
$\H_\gamma (H_{\mu^*}[X;q]) = H_{(\gamma,\mu^{(1)},
\ldots, \mu^{(k)})}[X;q]$ (see \cite{SZ})
\end{itemize}

In addition, analogs of most properties of the
Hall-Littlewood functions and the Kostka-Foulkes
polynomials also seem to hold (see for instance
\cite{SW}).

We should also mention that there is an analog of
several other formulas for the Hall-Littlewood functions.
It follows from the definition of the functions
$H_{\mu^*}[X;q]$
that
\begin{equation}
H_{\mu^*}[X;q] = \Omega[Z_n X]
\prod_{1 \leq i < j \leq n} \left(1- z_j/z_i\right)
\prod_{(i,j) \in Roots_\eta}
\frac{1}{1-q z_j/z_i}\coeff_{z^{\bar {\mu^*}}}.
\end{equation}
For $k>0$, if we define the operation,
\begin{equation}
\H(Z^k) P[X]
= P[X - (1-q) Z^*]
\Omega[Z X] \prod_{1 \leq i < j \leq k}
\left(1-z_j/z_i\right),
\end{equation}
where $Z^* = \sum_{i=1}^k \frac{1}{z_i}$, then
\begin{equation}
\H(Z^{\eta_1}) \H(Z^{\eta_2}) \cdots
\H(Z^{\eta_{\ell(\eta)}}) 1 = \Omega[Z_n X]
\prod_{1 \leq i < j \leq n} \left(1- z_j/z_i\right)
\prod_{(i,j) \in Roots_\eta}
\frac{1}{1-q z_j/z_i}
\end{equation}
and therefore $\H(Z^k) H_{\mu^*}[X;q] \coeff_{z^\gamma} =
H_{(\gamma,\mu^*)}[X;q]$.

This construction exists in complete analogy within the
$Q$-function algebra.  We will create a family of functions
in $\Gamma$ which are indexed by a sequence of strict
partitions.  Let $\mu^\ast = (\mu^{(1)}, \mu^{(2)}, \ldots,
\mu^{(k)})$ where $\mu^{(i)}$ is a strict partition and set
$\eta = (\ell(\mu^{(1)}), \ell(\mu^{(2)}), \ldots,
\ell(\mu^{(k)}))$. Define $Roots_\eta$ as before and then
define the function
\begin{equation}\label{Gmuparab}
G_{\mu^\ast}[X;q] = \prod_{(i,j) \in Roots_\eta}
\frac{1+q R_{ij}}{1-q R_{ij}} Q_{\bar {\mu^*}}[X].
\end{equation}

We may also view these elements of $\Gamma$ as
the result of a family of operators acting on $1$.
Consider the composition of the operators
\begin{align*}\QS_\la P[X]
&:= \QS_{\la_1} \QS_{\la_2} \cdots \QS_{\la_k} P[X]\\
&= P[X- Z_k^*] \Omega[(1-\minus)Z_k X]
\prod_{1 \leq i < j \leq k} \frac{1 - z_j/z_i}{1 + z_j/z_i}
\coeff_{z^\la}
\end{align*}
and then set $\G_\la := {\widetilde {\QS_\la}}^q$, that is
\begin{equation}
\G_\la P[X]
= P[X+(q-1) Z_k^*] \Omega[(1-\minus)Z_k X]
\prod_{1 \leq i < j \leq k} \frac{1 - z_j/z_i}{1 + z_j/z_i}
\coeff_{z^\la}.
\end{equation}
As with the other operators of this sort, it is easily
shown that a composition of $G_{\mu^{(i)}}$ acting on
$1$ is equivalent to the defining relation (\ref{Gmuparab})
of the functions $G_{\mu^*}[X;q]$ and hence that
\begin{equation}
\G_\gamma G_{\mu^*}[X;q] = G_{(\gamma, \mu^{(1)}, \ldots,
\mu^{(k)})}[X;q].
\end{equation}

The $G_{\mu^\ast}[X;q]$ functions seem to share many of the
same properties of the $H_{\nu^\ast}[X;q]$ and
$G_\gamma[X;q]$ analogs.
Define the polynomials
$L_{\la;\mu^\ast}(q)$ by the
expansion
\begin{equation}
G_{\mu^\ast}[X;q] = \sum_{\la} L_{\la;\mu^\ast}(q)
Q_\lambda[X].
\end{equation}

\begin{itemize}
\item
$G_{\mu^\ast}[X;0] = Q_{\bar {\mu^\ast}}[X]$.  $\bar
{\mu^\ast}$ need not be a strict partition, however if
it is not then the straightening relations (\ref{m-n}),
(\ref{m-m}) and (\ref{0}) may be applied to reduce
the expression.
\item
$G_{\mu^\ast}[X;1] = Q_{\mu^{(1)}}[X] Q_{\mu^{(2)}}[X]
\cdots Q_{\mu^{(k)}}[X]$ and hence $L_{\la;\mu^*}(1)$ is
equal to the coefficient of $Q_\la[X]$ in
the product $Q_{\mu^{(1)}}[X] Q_{\mu^{(2)}}[X]
\cdots Q_{\mu^{(k)}}[X]$.
\item If $\mu^* =
((\gamma_1),(\gamma_2), \ldots, (\gamma_{\ell(\gamma)}))$
where $\gamma$ is a partition, then
$G_{\mu^\ast}[X;q] = G_\gamma[X;q]$.
\item If $\bar {\mu^*}$ is a strict partition then
$G_{\mu^*}[X;q] = Q_{\bar {\mu^\ast}}[X] + \sum_{\la > {\bar
{\mu^*}}} L_{\la;\mu^*}(q) Q_\la[X]$.
\end{itemize}

Computing these coefficients suggests the following
remarkable conjecture and indicates that these coefficients
are an important
$q$-analog of the structure coefficients of the $Q_\la[X]$
functions in the same way that the $K_{\la;\mu^*}(q)$
polynomials are $q$-analogs of the Littlewood-Richardson
coefficients.

\begin{conj} \label{conjparab} For a sequence of partitions
$\mu^*$, if ${\bar {\mu^*}}$ is a partition
then $L_{\la;\mu^\ast}(q)$ is a
polynomial in $q$ with non-negative integer coefficients.
\end{conj}

If this conjecture is true then the polynomials
$L_{\la;\mu^*}(q)$ are a $q$ analog of the coefficient
of $Q_\la[X]$ in the product $Q_{\mu^{(1)}}[X]
Q_{\mu^{(2)}}[X]
\cdots Q_{\mu^{(k)}}[X]$. A combinatorial description
for these coefficients was given in \cite{S1} and hence
we are looking for an additional statistic on the set
of objects counted by them which includes as a special
case the coefficients $L_{\la\mu}(q)$.

This conjecture suggests that the
$L_{\la\mu^\ast}(q)$ should also share many of the
properties that are held by the $K_{\la;\mu^*}(q)$ and that
generalize the case of the Kostka-Foulkes polynomials.

We remark that the parabolic Kostka polynomials indexed
by a sequence of partitions $\mu^*$ where each $\mu^{(i)}$
is a rectangle (i.e. each $\mu^{(i)} = (a_i, a_i,\ldots, a_i)$ for
some $a_i$) is a special subfamily of these  polynomials.
In this case, explicit combinatorial formulas are known for
the coefficients (see for example \cite{SiWa}, \cite{Sh} or
\cite{KS}) which imply that the coefficients
$K_{\la;\mu^*}(q)$ are positive.  By contrast, for the
generalizations of the $Q$-Kostka polynomials we know that
if $\mu^{(i)}$ has two equal parts for any $i$ then
$G_{\mu^*}[X;q] = 0$, hence this special case is not of interest
in this setting.


\section{Appendix: Tables of $2^{\ell(\la)-\ell(\mu)} L_{\la\mu}(q)$ for 
$n=4,5,6,7,8,9$}

\begin{equation*}
\left [\begin {array}{cc} (3,1)&(4)\\1&q\\\noalign{\medskip}0&1\end 
{array}
\right ]
\end{equation*}

\begin{equation*}
\left [\begin {array}{ccc} 
(3,2)&(4,1)&(5)\\1&2\,q&{q}^{2}\\\noalign{\medskip}0&1&q
\\\noalign{\medskip}0&0&1\end {array}\right ]
\end{equation*}

\begin{equation*}
\left [\begin {array}{cccc} (3,2,1)&(4,2)&(5,1)&(6)\\
1&2\,{q}^{2}+q&2\,{q}^{3}+{q}^{2}&{q}^{4}
\\\noalign{\medskip}0&1&2\,q&{q}^{2}\\\noalign{\medskip}0&0&1&q
\\\noalign{\medskip}0&0&0&1\end {array}\right ]
\end{equation*}

\begin{equation*}
\left [\begin {array}{ccccc} (4,2,1)&(4,3)&(5,2)&(6,1)
&(7)\\1&q&2\,{q}^{2}+q&2\,{q}^{3}+{q}^{2}&{q}^{
4}\\\noalign{\medskip}0&1&2\,q&2\,{q}^{2}&{q}^{3}\\\noalign{\medskip}0
&0&1&2\,q&{q}^{2}\\\noalign{\medskip}0&0&0&1&q\\\noalign{\medskip}0&0&0
&0&1\end {array}\right ]
\end{equation*}

\begin{equation*}
\left [\begin {array}{cccccc} (4,3,1)&(5,2,1)&(5,3)&(6,
2)&(7,1)&(8)\\
1&2\,q&2\,{q}^{2}+q&2\,{q}^{2}+2\,{q}^{3
}&{q}^{3}+2\,{q}^{4}&{q}^{5}\\\noalign{\medskip}0&1&q&2\,{q}^{2}+q&2\,
{q}^{3}+{q}^{2}&{q}^{4}\\\noalign{\medskip}0&0&1&2\,q&2\,{q}^{2}&{q}^{
3}\\\noalign{\medskip}0&0&0&1&2\,q&{q}^{2}\\\noalign{\medskip}0&0&0&0&
1&q\\\noalign{\medskip}0&0&0&0&0&1\end {array}\right ]
\end{equation*}

\begin{equation*}
\left [\begin {array}{cccccccc} (4,3,2)&(5,3,1)&(5,4)&(6,2,1)&
(6,3)&(7,2)&(8,1)&(9)\\
1&2\,q+4\,{q}^{2}&2\,{q}^{3}+{q}^{2}&2
\,{q}^{2}+4\,{q}^{3}&{q}^{2}+2\,{q}^{4}+4\,{q}^{3}&4\,{q}^{4}+{q}^{3}+
2\,{q}^{5}&2\,{q}^{6}+2\,{q}^{5}&{q}^{7}\\\noalign{\medskip}0&1&q&2\,q
&2\,{q}^{2}+q&2\,{q}^{2}+2\,{q}^{3}&{q}^{3}+2\,{q}^{4}&{q}^{5}
\\\noalign{\medskip}0&0&1&0&2\,q&2\,{q}^{2}&2\,{q}^{3}&{q}^{4}
\\\noalign{\medskip}0&0&0&1&q&2\,{q}^{2}+q&2\,{q}^{3}+{q}^{2}&{q}^{4}
\\\noalign{\medskip}0&0&0&0&1&2\,q&2\,{q}^{2}&{q}^{3}
\\\noalign{\medskip}0&0&0&0&0&1&2\,q&{q}^{2}\\\noalign{\medskip}0&0&0&0
&0&0&1&q\\\noalign{\medskip}0&0&0&0&0&0&0&1\end {array}\right ]
\end{equation*}

\vskip .5in

\vfill
\pagebreak

\section{Appendix: example of conjectured tableaux 
poset of content $(4,3,2)$}

$$\Young{\Blk&\Blk&3&3\cr\Blk&2&2&2\cr1&1&1&1\cr}
$$

$$\Young{\Blk&\Blk&3\cr\Blk&2&2&2\cr1&1&1&1&3^*\cr}
$$

$$\Young{\Blk&\Blk&3\cr\Blk&2&2&3^*\cr1&1&1&1&2^*\cr}\hskip
.3in
\Young{\Blk&2&2&2&3^*\cr1&1&1&1&3'\cr}\hskip .3in
\Young{\Blk&\Blk&3\cr\Blk&2&3'\cr1&1&1&1&2^*&2\cr}\hskip
.3in
\Young{\Blk&2&2&2\cr1&1&1&1&3^*&3\cr}$$

$$\Young{\Blk&3&3\cr1&1&1&1&2^*&2&2\cr}
\hskip .3in
\Young{\Blk&2&2&3^*&3\cr1&1&1&1&2^*\cr}\hskip .3in
\Young{\Blk&\Blk&3\cr\Blk&2&2\cr1&1&1&1&2^*&3^*\cr}\hskip
.3in
\Young{\Blk&2&2&3^*\cr1&1&1&1&2^*&3^*\cr}
$$

$$\Young{\Blk&2&3^*\cr1&1&1&1&2^*&2&3^*\cr}
\hskip .3in
\Young{\Blk&2&3^*&3\cr1&1&1&1&2^*&2\cr}
$$

$$\Young{\Blk&2&2\cr1&1&1&1&2^*&3^*&3\cr}
\hskip .3in
\Young{\Blk&3\cr1&1&1&1&2^*&2&2&3^*\cr}
$$

$$\Young{\Blk&2\cr1&1&1&1&2^*&2&3&3^*\cr}$$

$$\Young{1&1&1&1&2^*&2&2&3^*&3\cr}
$$

\font\Sc=cmcsc10

\noindent
{\Sc Figure 2.} 
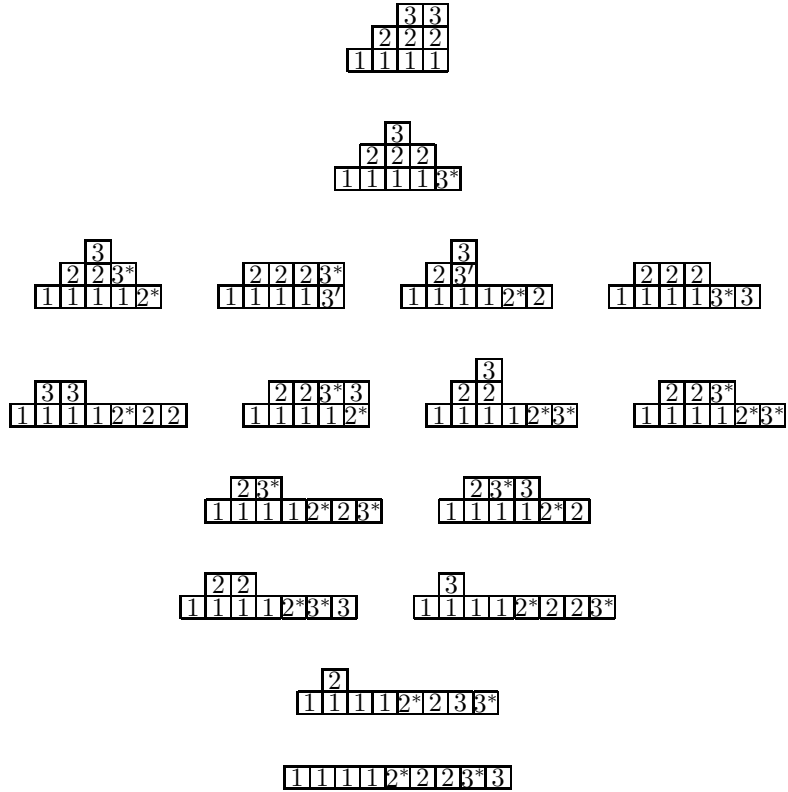
{The cells marked with a $k^*$ can be labeled with
either $k$ or $k'$, we conjecture that the statistic is independent of these
markings.  The value of
$G_{(4,3,2)}[X;q]$ determines the position of each of the
shifted tableaux here except for the two of shape $(8,1)$,
however the statistics in smaller polynomials
(e.g. $G_{(4,3,1)}[X;q]$) suggest this rank function.  
The covering relation is unknown, but the rank function
indicates that it is not the same as the charge statistic.}

\vskip .1in
\noindent
{\Sc Acknowledgment:} Thank you to Nantel Bergeron for many helpful
suggestions on this research.

\end{document}